\documentclass[12pt]{article}
\usepackage{amsmath,amsthm,amsfonts,amssymb}
\usepackage[all]{xy}
\usepackage[usenames]{color} 
\usepackage{pdfsync}

\usepackage{rotating}

\newcounter{enumitemp}

\newcommand\pref[1]{(\ref{#1})}

\newtheorem{thm}{Theorem}[section]

\newtheorem{theorem}[thm]{Theorem}

\newtheorem{lemma}[thm]{Lemma}

\newtheorem{cor}[thm]{Corollary}
\newtheorem{corollary}[thm]{Corollary}
\newtheorem{proposition}[thm]{Proposition}
\newtheorem*{proposition*}{Proposition}
\newtheorem{question}[thm]{Question}

\newtheorem{fact}[thm]{Fact}

\theoremstyle{definition}

\newtheorem{definition}[thm]{Definition} 

\newtheorem*{defn*}{Definition}

\newtheorem{notation}[thm]{Notation}

\newtheorem{remark}[thm]{Remark}
\theoremstyle{remark}

\newcounter{remarks}
{\paragraph*{Remarks}\smallskip
 \begin{list}{\arabic{remarks}. }{\usecounter{remarks}%
 \setlength{\leftmargin}{0in}%
 \setlength{\rightmargin}{0in}%
 \setlength{\labelsep}{0pt}%
 \setlength{\labelwidth}{0pt}%
 \setlength{\listparindent}{0pt}%
 }
}
{
\end{list}
}

\newcommand\from\colon

\newcommand\subgroup{<}
\newcommand\normal\triangleleft
\newcommand\infinity\infty

\newcommand\supp{\text{supp}}
\newcommand\disjunion\coprod
\newcommand\act\curvearrowright

\DeclareMathOperator{\Fix}{Fix}
\DeclareMathOperator{\Per}{Per}

\newcommand{\R}{{\mathbb R}}

\newcommand{\Z}{{\mathbb Z}}

\newcommand{\E}{{\mathcal E}}
\newcommand{\V}{{\mathcal V}}

\newcommand{\K}{{\mathcal K}}

\renewcommand\P{{\mathcal P}}
\newcommand{\A}{\mathcal A}

\renewcommand\L{\mathcal L}
\newcommand{\PF}{{\text{PF}}}

\DeclareMathOperator{\Out}{\mathsf{Out}}
\DeclareMathOperator{\Aut}{\mathsf{Aut}}

\DeclareMathOperator{\Stab}{\mathsf{Stab}}

\newcommand{\pg}{PG}

\newcommand{\upg}{UPG}
\newcommand{\F}{\mathcal F}
\newcommand{\rtt}{relative train track map}

\newcommand{\fG} {f : G \to G}
\newcommand{\ti} {\tilde}
\newcommand{\iNp} {indivisible Nielsen path}

\newcommand{\noneg}{NEG}
\renewcommand\neg\noneg

\newcommand{\wt}{\widetilde}

\newcommand{\ct}{CT}
\newcommand{\cts}{CTs}

\newcommand{\comment}[1]{}

\newcommand\BookOneTag{BFH:TitsOne}
\newcommand\BookOne{\cite{\BookOneTag}}
\newcommand\BookTwoTag{BFH:TitsTwo}
\newcommand\BookTwo{\cite{\BookTwoTag}}
\newcommand\BookThree{\cite{BFH:Solvable}}

\newcommand\recognitionTag{FeighnHandel:recognition}
\newcommand\recognition{\cite{\recognitionTag}}
\newcommand\abelianTag{FeighnHandel:abelian}
\newcommand\abelian{\cite{\abelianTag}}

\newcommand\SubgroupsOneTag{HandelMosher:SubgroupsI}

\newcommand\SubgroupsTwoTag{HandelMosher:SubgroupsII}
\newcommand\SubgroupsTwo{\cite{\SubgroupsTwoTag}}

\newcommand\bdy\partial

\newcommand\intersect\cap
\newcommand\union\cup
\newcommand\<\langle
\renewcommand\>\rangle
\newcommand\meet\wedge

\newcommand\cross\times
\newcommand\restrict{\bigm |}
\newcommand\wh{\widehat}
\newcommand\inject\hookrightarrow

\newcommand\injectto\hookrightarrow
\newcommand\injectfrom\hookleftarrow
\newcommand\surjectto\twoheadrightarrow
\newcommand\surjectfrom\twoheadleftarrow

 \newcommand\surjection\twoheadrightarrow
 
\DeclareMathOperator\gl{GL}

\DeclareMathOperator\IA{IA}

\newcommand\IAThree{\IA_n(\Z/3)}

\newcommand\cH{{\cal H}}

\title{Virtually abelian subgroups of $\IAThree$ are abelian}

\author{Michael Handel\thanks{The first author was supported by 
  National Science Foundation grant  DMS-1308710 and by   a PSC-CUNY  grant in Program Year 47.} and Lee Mosher\thanks{The second author is supported by National Science Foundation grant DMS-1708361.}}

\begin{document}
\maketitle
\begin{abstract}  When studying subgroups of $\Out(F_n)$,  one often  replaces a given subgroup $\cH$ with one of its finite index subgroups $\cH_0$  so that virtual properties of $\cH$ become actual properties of $\cH_0$.  In many cases,   the finite index subgroup is  $\cH_0 = \cH \cap \IAThree$.     For which properties is this a good choice?  
  Our main theorem states that being  abelian is such a property.  Namely, every virtually abelian subgroup of $\IAThree$ is abelian.
\end{abstract}

\tableofcontents

\section {Introduction}

It is common, when studying elements  of  $  \Out(F_n)$, to replace the given element by an iterate in order  to improve its invariance properties.     For example,  each $\theta \in \Out(F_n)$  has an iterate $\phi = \theta^k$ satisfying the following properties.
\begin{enumerate}
\item If some iterate of $\phi$ fixes a conjugacy class $[a]$ then $\phi$ fixes $[a]$.
\item If some iterate of $\phi$ fixes the conjugacy class [F] of a free factor $F$ then $\phi$ fixes $[F]$.
\item $\phi$ fixes   each element in its set $\L(\phi)$ of   attracting laminations.
\item  $\phi$ fixes   each element in its set of singular rays and eigenrays.
\end{enumerate}
 If $\theta$ is {\em rotationless} in the sense of \cite{\recognitionTag} then iteration is not necessary: each of the above properties is automatically satisfied by $\phi = \theta$   \cite[Lemma 3.30 and Definition 3.13]{\recognitionTag}.    Every $\theta$ has a rotationless iterate and the number of iterates required is uniformly bounded \cite[Lemma 4.42]{\recognitionTag}.

The subgroup analog of replacing an individual element $\theta$  with a rotationless iterate $\theta^k$ is to replace a given subgroup $\cH$ with its finite index subgroup $\cH \cap \IAThree$ where    $\IAThree\subgroup \Out(F_n)$ is the finite index subgroup  consisting of elements that act trivially on $\Z/3$ - homology.    This was done, for example, in the proof of   the Tits Alternative for $\Out(F_n)$ \BookOne, \BookTwo, \BookThree.

  In \SubgroupsTwo\ 
    (see also \cite[Propositions 3.16 and  4.41]{\BookTwoTag})  we proved       that all elements of $\IAThree$ 
    satisfy (1) - (3) above.   (If an element of $\IAThree$ satisfies (4)  then it is rotationless \cite[Lemma 3.12]{FeighnHandel:ctAlgorithm}.)  These invariance properties played a significant role in our    series of papers \cite{HandelMosher:Subgroups}  establishing the \lq subgroup decomposition\rq\ theorem for $\Out(F_n)$ and then again  in  (\cite{HandelMosher:BddCohomologyI}, \cite{HandelMosher:BddCohomologyII}), in which the  $H^2_b$-alternative for $\Out(F_n)$  is established: for every  finitely generated subgroup  $\cH \subgroup \Out(F_n)$ either $\cH$ is virtually abelian or  $H^2_b(\cH;\R)$ has uncountably infinite dimension.  

 The main result of this paper is motivated in part by  \cite{HandelMosher:BddCohomologyI} and\cite{HandelMosher:BddCohomologyII}, in which virtually abelian subgroups appear naturally and in which information is lost when one passes to finite index subgroups, and in part by our appreciation of the importance of   $\IAThree$.  Having seen that   elements of $\IAThree$ satisfy (1) - (3) without iteration, one can ask analogously, which virtual properties of arbitrary subgroups of $\Out(F_n)$ are true for subgroups of $\IAThree$ without passing to a subgroup of finite index?   Our main theorem in this paper is one such property.

\begin{theorem} \label{main} Each  virtually abelian subgroup $\cH \subgroup \IAThree$ is abelian.
\end{theorem}

Abelian subgroups of $\Out(F_n)$ are finitely generated and    $\IAThree$ is torsion free; the former is contained in   \cite{BassLubotzky:Linear-central}   and the latter follows  from  \cite[Corollary 5.7.6]{\BookOneTag}.   Thus,

\begin{corollary}  Every virtually abelian subgroup of $\Out(F_n)$ has a finitely generated,   free abelian subgroup of index  at most $|GL(n, \Z_3)| < 3^{n^2}$. 
\end{corollary}

       In Section~\ref{reduction},   after a brief review of PG and UPG subgroups, we reduce Theorem~\ref{main} to the following proposition.    

\begin{proposition}  \label{upg prop} Suppose that $K \subgroup \IAThree$ is an  abelian \upg\ subgroup.  Then the normalizer of $K$ in $\IAThree$ equals the centralizer of $K$ in $\IAThree$. 
\end{proposition}
All \pg\ elements of $\IAThree$\ are \upg\     \cite[Corollary 5.7.6]{\BookOneTag} and, in fact, rotationless (Lemma~\ref{upg is rotationless}). Therefore Proposition~\ref{upg prop} may be equivalently restated using \pg\ in place of \upg.
The proof of Proposition~\ref{upg prop} appears in Section~\ref{main proof}.

\medspace
Continuing with the theme of studying $\IAThree$, we pose the following question, the answer to which is yes if $\phi$ and $\psi$ are rotationless by an easy application of \cite[Theorem 5.3]{\recognitionTag}.

\begin{question} Are roots unique in the group $\IAThree$?  That is, if $\phi, \psi \in \IAThree$ and $\phi^k =\psi^k$ for some $k\ge1$, is $\phi = \psi$?
\end{question}

  Section~\ref{sec:background} contains background material including subsections on    \upg\ elements and \upg\ subgroups. 
  
\section{Reduction to  Proposition~\ref{upg prop}} \label{reduction}

 Each $\psi \in \Out(F_n)$ has an associated finite set   $\L(\psi)$ of attracting laminations, each of which is invariant under some iterate of $\psi$ \cite[Section 3.1]{\BookOneTag}.  For a subgroup $\cH \subgroup \Out(F_n)$, we let $\L(\cH)=  \cup_{\psi \in \cH} \L(\cH)$.    
  If   $\L(\psi) = \emptyset$, then we say that $\psi$ has {\em polynomial growth} and write $\psi \in \pg(F_n)$ or simply $\psi \in \pg$.   If in addition, the image of  $\psi$  in $\gl(n,\Z)$ is  unipotent  
     then we write $\psi \in \upg(F_n)$ or simply $\theta \in \upg$ \BookOne ,   \BookTwo.      
     
\medspace

\noindent{\em Proof of Theorem~\ref{main} assuming Proposition~\ref{upg prop}:}\ \  Let $\cH \subgroup \IAThree$ be virtually abelian.  We first follow the proof of \cite[Theorem 7.0.1]{\BookOneTag} to show that that there is an exact sequence 
$$ 1 \to K \to \cH \to \Z^k \to 1$$
for some $k$  and some abelian subgroup $K \subgroup \upg$.    By \cite[Lemma 4.7]{HandelMosher:BddCohomologyII},  $\L(\cH)$ is a finite collection $  \{\Lambda_1,\ldots, \Lambda_k\}$ of $\cH$-invariant laminations.   For each $1 \le i \le k$, let $\PF_{\Lambda_i} : \Stab(\Lambda_i) \to \Z$ be the expansion  factor homomorphism \ for $\Lambda_i$ as defined in  \cite[Section 3.3]{\BookOneTag}.   Let   $\PF = \oplus_{i=1}^k\PF_{\Lambda_i} : \cH \to \Z^k$   be the direct sum of the restrictions to $\cH$ of the $\PF_{\Lambda_i}$'s.        If $\theta$ is an element of the kernel $K$ of $\PF$  then $\L(\theta) \cap  \L(\cH) = \emptyset$ by \cite[Corollary 3.3.1]{\BookOneTag}. Thus  $\L(\theta) = \emptyset$    and  $K$ is \pg.    Applying our assumption that $\cH \subgroup \IAThree$, we have that  $K$ is \upg\ by  \cite[Corollary 5.7.6]{\BookOneTag}.   It then follows  that $K$ is solvable \cite[Corollary 1.3]{\BookTwoTag}.   Since $K$ is virtually abelian, it is finitely generated by \cite{BassLubotzky:Linear-central}  (see also \cite{BFH:Solvable}).  
We  can therefore apply \cite[Corollary 3.11]{BFH:Solvable} to conclude that $K$ is abelian.

  Proposition~\ref{upg prop} implies that $K$ is in the center of $\cH$.    In particular $ [\psi_1,\psi_2] $  (which is an element of $K$)  commutes with  $\psi_1$ and $\psi_2$   for all $\psi_1, \psi_2\in \cH$.   For all $p \ge 1$ and all $\psi_1,\psi_2 \in \cH$ we have  $[\psi_1,\psi_2]^p = [\psi_1^p,\psi_2]$ and  similarly  $[\psi_1,\psi_2]^p =  [\psi_1,\psi_2^p]$.  For the first of these equations the inductive step is:  
  \begin{align*}
  [\psi_1,\psi_2]^p = [\psi_1,\psi_2]^{p-1}  \psi_1 \psi_2 \psi_1^{-1} \psi_2^{-1} & =  \psi_1[\psi_1,\psi_2]^{p-1}   \psi_2 \psi_1^{-1} \psi_2^{-1} \\&= \psi_1[\psi_1^{p-1},\psi_2]    \psi_2 \psi_1^{-1} \psi_2^{-1} \\ & =\psi_1\psi_1^{p-1}\psi_2\psi_1^{1-p} \psi_2^{-1} \psi_2 \psi_1^{-1} \psi_2^{-1}   = [\psi_1^p,\psi_2]
  \end{align*}
   Since $\cH$ is virtually abelian, there exists $p\ge 1$ such that $ [\psi_1^p,\psi_2^p]$ is trivial.  It follows that $[\psi_1,\psi_2]^{p^2}$ is trivial. Since finite order \upg\ elements are trivial \cite[Lemma 4.47]{\BookTwoTag}, we conclude that    $[\psi_1,\psi_2]$ is trivial for all $\psi_1, \psi_2 \in \cH$.  
\qed

\section{Background} \label{sec:background}
\subsection{Basics}

Much of the material in this subsection is standard and is included   to establish notation and for convenient reference.  Further details can be found in \cite[Section~2]{\BookOneTag}, \cite[Section~2]{\recognitionTag} or   \cite[Section 1]{\SubgroupsOneTag}. 

\paragraph{Marked graphs} The free group $F_n$ of rank $n$ is identified with $\pi_1(R_n)$ where $R_n$ is the graph with one vertex and $n$ edges.  A {\em marked $n$-graph} is a connected finite graph $G$ of rank $n$ that has no valence one vertices and is  equipped with a homotopy equivalence $ R_n \to G$ called a {\em marking of $G$}.  The marking provides an identification of $F_n$ with $\pi_1(G)$ that is well defined up to inner automorphism.  A homotopy equivalence $\fG$ determines an outer automoprphism of $\pi_1(G)$ and hence an element $\phi \in \Out(F_n)$ that we say is {\em represented by $\fG$}.   

 Edges of $G$ are assumed to be oriented with $\bar E = E^{-1}$ denoting the edge $E$ with its orientation reversed. All of the $\fG$ that we consider will take vertices to vertices and restrict to an immersion on each edge.   
A {\em direction} $d$ at a vertex $v \in G$ is the germ of an oriented edge $E$  with initial vertex $v$.      Define the action of $f$ on directions by $d \mapsto d'$ where $d'$ is the  direction determined by the first edge in $f(E)$.   

We denote the universal cover of $G$   by $\ti G$ and the set of ends of $\ti G$ by $\partial \ti G$.  
\begin{fact}  \label{f hat} Suppose that $\fG$ is a homotopy equivalence of a marked graph.  Then each lift $\ti f:\ti G \to \ti G$ extends continuously over $\partial \ti G$ by a homeomorphism $\hat f :\partial \ti G \to \partial \ti G$.
\end{fact}  

For each $a \in F_n$, the inner automorphism $i_a :F_n \to F_n$ is defined by $i_a(g) = aga^{-1}$;  the conjugacy class of $a$ is denoted by $[a]$.

\begin{fact} \label{theta hat} Each $\Theta \in\Aut(F_n)$  extends continuously  to  a homeomorphism $\wh \Theta: \partial F_n \to \partial F_n$. For each non-trivial inner automorphism $i_a$,  its boundary extension $\wh{i_a}$  fixes two points, {\em a source $a^-$} and {\em a sink $a^+$}. 
\end{fact}

\begin{fact}\label{boundary identification} For each marked graph $G$, the identification of $\pi_1(G)$ with $F_n$ induces 
\begin{enumerate}
\item \label{marking1} an identification of the group of covering translations of $\ti G$ with $F_n$ and
\item \label{marking2} an identification of  $\partial \ti G$ with $\partial F_n$ 
\end{enumerate} 
so that  for each non-trivial $a \in F_n$, if $T_a$ is the covering translation of $\ti G$ identified by \pref{marking1} with  $a$, and if $A_a$ is the axis of $T_a$ then the following hold:    $a^-$ and $a^+$ are   identified by \pref{marking2} with  the repelling and attracting  endpoints  of  $A_a$  respectively;  and the projection of $A_a$ into $G$ is a {\em circuit that represents the conjugacy class $[a]$ of $a$}.   More generally, for any $\theta\in \Out(F_n)$ and homotopy equivalence $\fG$ representing  $\theta$, there is  
 a bijection $$\ti f \leftrightarrow \Theta$$ between the set of lifts $\ti f : \ti G\to \ti G$  and the set of $\Theta \in \Aut(F_n)$ representing $\theta$ such that the homeomorphisms $\hat f : \partial \ti G \to \partial \ti G$ and $\wh \Theta : \partial F_n \to \partial F_n$   defined in Fact~\ref{f hat} and Fact~\ref{theta hat} agree under  identification  \pref{marking2}.
\end{fact}

 \begin{fact} \label{fix both} \cite[Lemma~2.4] {BFH:Solvable}  Suppose that $\Phi \in \Aut(F_n)$ and $a\in F_n$. If $\partial \Phi$ fixes either $a^+$ or $a^-$ then $a \in \Fix(\Phi)$ and $\partial \Phi$  fixes both $a^+$ and $a^-$.
 \end{fact}

\paragraph{Principal lifts, rotationless outer automorphisms and rotationless maps \cite{\recognitionTag}} We will only be interested in principal lifts and principal vertices in the \upg\ setting and so we can give simplified versions of their definitions.

\begin{fact}  \cite[Proposition I.1]{GJLL:index} \label{gjll}
For  $\Theta \in \Aut(F_n)$, we denote the  fixed subgroup of $\Theta$ by  $\Fix(\Theta)$.    
   There is a disjoint union  
$$\Fix(\wh \Theta) = \Fix_-(\wh \Theta) \cup \Fix_+(\wh \Theta) \cup \partial \Fix(\Theta)$$
 where $\Fix_-(\wh \Theta) \subset \partial F_n$ is a finite union of $\Fix(\Theta)$-orbits of repellers and $\Fix_+(\wh \Theta) \subset \partial F_n$ is a finite union of $ \Fix(\Theta)$-orbits of   attractors. $\Fix_N(\Theta) \subset \partial F_n$ is defined to be $\partial \Fix(\Theta) \cup \Fix_+(\Theta)$.
 \end{fact}
 
 \begin{remark}  In the special case that $\Fix(\Theta) = \langle a \rangle$ and $\Fix(\wh \Theta) = \partial \Fix(\Theta) = \{a^{\pm}\}$, it may happen that $a^+$ or $a^-$  has an attracting neighborhood   for the action of $\wh \Theta$.  This happens for example if $\Theta = i_a$.  In all other cases,   $\Fix_+(\wh \Theta)$ is exactly the set of isolated attractors and $\Fix_-(\wh \Theta)$ is exactly the set of isolated repellers.  
 \end{remark}
 
 \begin{notation}\label{notn:eigenray} The  $F_n$-orbit of an element of $\Fix_+(\wh \Theta)$ is  called an {\em eigenray} for $\theta$.   The $F_n$-orbit of an element of $\Fix_-(\wh \Theta)$  is an  eigenray for $\theta^{-1}$.
\end{notation}
\begin{definition} \label{def:principal}  An automorphism $\Theta  $ representing a \upg\  $\theta \in \Out(F_n)$  is {\em principal}  \cite[Definition 3.1]{\recognitionTag} if either $\Fix_N(\Theta)$ contains at least three points or if $\Fix_N(\Theta)$ is a two point set that is not   $\{a^\pm\}$ for  some $a \in F_n$ on $\partial F_n$.  The set of principal $\Phi$ representing $\theta$ is denoted by $P(\theta)$.  An element $\phi \in \Out(F_n)$ is {\em rotationless} \cite[Definition 3.13]{\recognitionTag} if: (i) $\Phi \mapsto \Phi^k$  defines a bijection between $P(\phi)$ and $P(\phi^k)$ for all $k \ge 1$; and (ii) $\Fix_N(\Phi) = \Fix_N(\Phi^k)$ for all $\Phi \in P(\phi)$ and all $k \ge 1$.

 If $\fG$ represents $\theta$ and $\ti f  : \ti G \to \ti G$ corresponds to $\Theta \in P(\theta)$ as in Fact~\ref{boundary identification} then we say that {\em $\ti f$ is principal}.  An element $x$ of  the set  $\Per(f)$ of $f$-periodic points  is {\em principal}  \cite[Definition 3.18]{\recognitionTag} unless it is contained in a component $C$ of $\Per(f)$ that is topologically a circle and each point in $C$ has exactly two periodic directions.  If each principal vertex and periodic direction at a principal vertex has period one then we say that $\fG$ is {\em rotationless}.
\end{definition}

\paragraph{Paths}   An {\em edge path} $\sigma$ in a marked graph $G$   is a concatenation of edges $\sigma = \ldots E_i E_{i+1} \ldots$ of $G$  where the terminal endpoint of $E_i$ equals the initial endpoint of $E_{i+1}$ for all $i$.  If there is no backtracking, i.e.  if  $E_{i+1} \ne E_i^{-1}$ for all $i$,  then we say that  $\sigma$ is a {\em path}.  If a path $\sigma \subset G$ is a bi-infinite concatenation then we say that $\sigma$ is a {\em line} in $G$. (All of the lines in this paper are oriented.) If a path $\sigma \subset G$ is a singly infinite concatenation then we say that $\sigma$ is {\em ray}.  We also allow the {\em trivial path} which is just a single vertex.   Concatenation $\sigma = \sigma_1 \sigma_2$ of edge paths $\sigma_1$ and $\sigma_2$ is defined if $\sigma_1$ has a terminal vertex, $\sigma_2$ has an initial vertex and if these vertices are equal.  The concatenation of paths need not be path.   

Paths and edge paths in $\ti G$ are defined similarly. Edge paths in $G$ lift to edge paths in $\ti G$ with paths lifting to paths.  A line in $G$ lifts to a line in $\ti G$ with well defined distinct ideal endpoints in $\partial \ti G$.  Conversely, every ordered pair of distinct points in $\partial \ti G$ is the ideal endpoint pair for a unique line in $\ti G$. A ray in $G$ lifts to a ray in $\ti G$ with one endpoint at a vertex and the other an ideal endpoint in $\partial \ti G$.   

Suppose that $f :G \to G$ is a homotopy equivalence  and    $\ti f : \ti G \to \ti G$ is a lift.  For any finite path  $\ti \sigma \subset \ti G$ with endpoints $\ti x,\ti y$,   we define $\ti f_\#(\ti \sigma)$ to be the unique path   with endpoints $\ti f(\ti x), \ti f(\ti y)$.   We define $\ti f_\#(\ti \sigma)$ for rays and lines similarly using $\hat f$ if one or both endpoint is ideal.    This descends to a well defined action  $\sigma \mapsto f_\#(\sigma)$ of $f$ on the set of  paths in $G$.   

 A {\em circuit} in $G$ is a cyclic concatenation of edges without backtracking and so can be viewed as an immersion of a circle.  A circuit in $G$ lifts to a line in $\ti G$ and we can extend the definition of $f_\#$ to include circuits.  A closed path $\sigma$ determines a circuit if the initial edges of $\sigma$ and $\bar \sigma$ are distinct.     If a circuit $\sigma$ represents the conjugacy class $[a]$ of $a \in F_n$ and if $\fG$ represents $\theta\in \Out(F_n)$ then $f_\#(\sigma)$ represents $\theta([a])$.

 A decomposition $\sigma = \ldots \sigma_i \sigma_{i+1} \ldots$ into  subpaths is a {\em splitting} if   $$f^k_\#(\sigma) =  \ldots f^k_\#(\sigma_i) f^k_\#(\sigma_{i+1})\ldots$$  is a decomposition into subpaths for all $k \ge 1$,      When the decomposition into $\sigma_i$'s is a splitting we write  $\sigma = \ldots \cdot \sigma_i \cdot\sigma_{i+1}\cdot  \ldots$

  An {\em abstract line} is the $F_n$-orbit of an ordered pair of distinct points in $\partial F_n$.  If $G$ is any marked graph then the identification of $\partial G$ with $\partial F_n$ (Fact~\ref{boundary identification}) defines a bijection between abstract lines and $F_n$-orbits of  lines in $\ti G$ and so also a bijection between abstract lines and lines in $G$.  An {\em abstract ray} is an $F_n$-orbit of a point in $\partial F_n$.  There is a bijection between abstract rays and equivalence classes of rays in $G$, where two rays in $G$ are equivalent if they have a common infinite subray.

\paragraph{Free factor systems}   \cite[Section 2.6]{\BookOneTag} If $A_1, \ldots,  A_k$ are non-trivial free factors and $A_1 \ast \ldots \ast A_k$ is a  free factor of $F_n$  then the set of conjugacy classes $ \{[A_1],\ldots, [A_k]\}$ is a   {\em free factor system}.    We write $$  \{[B_1],\ldots, [B_l]\} \sqsubset  \{[A_1],\ldots, [A_k]\}$$ and say that  $\{[B_1],\ldots, [B_l]\}$ is contained in $\{[A_1],\ldots, [A_k]\}$ if for each $B_i$ there exists $A_j$ so that some conjugate of $B_i$ is a subgroup of $A_j$.

   For every  inclusion $H \subset G$ of   a subgraph   in a  marked graph, there is an associated free factor system $\F(H) = \{[\pi_1(C_1)], \ldots, [\pi_1(C_k)]\}$ where  $\{C_1,\ldots,C_k\}$  is the set of  non-contractible components of $G$; see \cite[Example 2.6.1]{\BookOneTag} for details.  We say that  $H \subset G$ {\em realizes $\F(H)$}.   Every free factor system is realized by some $H \subset G$ and every nested sequence $\F_1 \sqsubset \F_2 \sqsubset \ldots \sqsubset \F_l$ is realized by some nested sequence of subgraphs $H_1 \subset H_2 \subset \ldots \subset H_l \subset G$.   One may assume without loss that the $H_i$'s are {\em core} subgraphs, meaning that all vertices have valence at least two. If $\F \sqsubset \F'$ can be realized by core subgraphs $H \subset H'$ such that $H' \setminus H$ is a single edge then we say that $\F \sqsubset \F'$ is {\em a one-edge extension};  otherwise, $\F \sqsubset \F'$ is a {\em multi-edge extension}.
   
\begin{fact} \cite[Section 2.6]{\BookOneTag} \label{restriction}  Suppose that $[F]$ is a $\theta$-invariant free factor conjugacy class and that $\Theta \in \Aut(F_n) $  represents $\theta$ and preserves $F$.   Then the element $\theta\restrict F$ of $\Out(F)$ determined by the restriction $\theta \restrict F$ is independent of the choice of $\Theta$.   
 \end{fact}  

A conjugacy class is {\em carried} by $[F]$ if some representative of it is an element of $F$.   An abstract ray   is {\em carried}  by $[F]$ if it is represented by a point in   $\partial F$.  An abstract line  is {\em carried}  by $[F]$ if it is represented by an ordered pair of points, both of which are contained   in  $\partial F$.  A conjugacy class, abstract ray or abstract line is carried by a free factor system $\F$ if it is carried by a component of $\F$.   If $H $  is a subgraph of a marked graph $G$ then a  conjugacy class [resp. abstract line] is carried by $\F(H)$ if and only if the corresponding circuit [resp. line] in $G$      is contained in $H$.  

\begin{fact} \cite[Fact 1.10]{\SubgroupsOneTag} (see also \cite[Section 2.6]{BestvinaHandel:tt})  \label{ffs support} For any set  $X$ of abstract lines, abstract rays and conjugacy classes  there is a unique minimal (with respect to $\sqsubset$) free factor system $\F_\supp(X)$ that carries each element of $X$.  If $\theta \in \Out(F_n)$ and $X$ is $\theta$-invariant then $\F_\supp(X)$ is $\theta$-invariant.
\end{fact}

\subsection{\upg\ elements}  

In this section we review some facts about individual \upg\ elements of $\Out(F_n)$. 

A \ct\ $\fG$ is a particularly nice kind of topological representative of $\theta \in \Out(F_n)$.   The complete definition of a \ct\ is given on \cite[page 47]{\recognitionTag}.  Since  we will only use \ct\  representatives in the special case when $\theta$ is  \upg, the definition can be simplified considerably.  Fact~\ref{ct basics} and the proof of Lemma~\ref{upg is rotationless} give a pretty complete picture of \ct s in this context.

 We delay the proof that every \upg\ element is rotationless, and hence represented by a \ct \ $\fG$ \cite[Theorem 4.28]{\recognitionTag}, until we have listed some properties enjoyed by such  \cts. 

 A \ct\ $\fG$ is equipped with a filtration $\emptyset = G_0 \subset G_1 \subset \ldots \subset G_N = G$ by $f$-invariant subgraphs. The subgraphs $H_r = G_r \setminus G_{r-1}$ are called the {\em strata}.   A path has {\em height $r$} if it is contained in $G_r$ and crosses at least one edge in $H_r$.   The set of fixed points for $f$  and the set of   periodic points for $f$ are denoted by $\Fix(f)$ and $\Per(f)$ respectively.      The set of vertices of $G$ is $V$.  Recall that $V$ is invariant under each $\fG$ that we consider.   .

\begin{fact} \label{ct basics}  Each \ct \ $\fG$ representing $\theta \in \upg$ satisfies the following properties.
\begin{enumerate}
\item \label{item:edge in ct} Each stratum $H_i$ is a single edge $E_i$. If $E_i$ is not fixed then there is a non-trivial closed path  $u_i \subset G_{i-1}$   such that $f(E_i) = E_i  u_i$.  
\item \label{item:fix} $ \Fix(f)= \Per(f) $ is the union of $V$ with the set of fixed edges.    
\item A direction based at a vertex is fixed if and only if it is periodic if and only if it is not  the terminal direction of a non-fixed edge.  
\item  Each vertex is principal.
\item \label{principal lifts} A lift $\ti f : \ti G \to \ti G$ is principal if and only if $\Fix(\ti f) \ne \emptyset$.    
\end{enumerate}
\end{fact}

\proof   The strata of $\fG$ are classified into three types: EG, NEG, and zero strata.  From \cite[Lemma 3.1.9]{\BookOneTag} and our assumption that $\L(\theta) = \emptyset$, it follows that $\fG$ has no EG strata.  The (Zero Strata) property of a \ct\ therefore implies that  $\fG$ has no zero strata.  Thus,   every stratum of $\fG$ is \noneg.     
Item (1) therefore follows from \cite[Lemma 4.21]{\recognitionTag}.  Items (2) and (3) follow from (1).  A vertex that is incident to  a fixed edge or is the terminal endpoint of a non-fixed edge is principal by the (Periodic Edges) and (Vertices) properties of a \ct\ respectively.  All other vertices are the initial endpoints of at least two non-fixed edges  and so are  principal by Definition~\ref{def:principal}.    This proves   (4).     Item  (5) follows from (4) and \cite[Remark 4.8, Corollaries 3.17 and 3.27]{\recognitionTag}.   
\endproof

A finite path $\sigma \subset G$ is a {\em Nielsen path} if $f_\#(\sigma) = \sigma$  and is an \iNp\ if there is no non-trivial decomposition of $\sigma$ into Nielsen subpaths.  Note that by Lemma~\ref{ct basics}\pref{item:fix}, we would have the same set of   \iNp s if we allowed paths to have endpoints that are not vertices.   

\medskip

   In order to apply  \ct\ theory to \upg\ elements we must prove that they are rotationless.  We will do this indirectly by using a result from \cite{\BookOneTag} to find   a pretty good \rtt, namely one that satisfies various of the conclusions of Fact~\ref{ct basics}, and then we will quote \cite[Proposition 3.29]{\recognitionTag}.   

\begin{lemma} \label{upg is rotationless}  Each $\theta \in \upg$ is rotationless.   
\end{lemma}
 
\proof    By \cite[Proposition 5.7.5]{\BookOneTag}, $\theta$ is represented by a \rtt\ $\fG$  and filtration $\emptyset = G_0 \subset G_1 \subset \ldots \subset G_N = G$  with a  subsequence of invariant core subgraphs $\emptyset = G_0 =G_{r(0)}  \subset G_{r(1)} \subset \ldots \subset G_{r(m)} =   G$ such that $G_{r(j+1)}$ is obtained from $G_{r(j)}$ in one of the following ways.    
\begin{description}
\item [(a)]  Adding a single fixed edge that is either a loop or has both endpoints in $G_{r(j)}$; $ r(j+1) = r(j) +1$.
\item [(b)] Adding a single  non-fixed  edge  satisfying Fact~\ref{ct basics}(1) with both endpoints in $G_{r(j)}$;  $r(j+1) = r(j) +1$.
\item [(c)] Adding two  non-fixed   edges satisfying Fact~\ref{ct basics}(1) with a common initial vertex not in $G_{r(j)}$ and both terminal endpoints in $G_{r(j)}$;  $r(j+1) = r(j) +2$.  
\end{description}

It is obvious from (a) - (c) that Fact~\ref{ct basics}(1)  is    satisfied.  This implies items (2) and (3) of Fact~\ref{ct basics}, which in turn prove that $\fG$ is rotationless  (Definition~\ref{def:principal}).    By \cite[Proposition 3.29]{\recognitionTag}, we are reduced to showing that $\fG$ satisfies the five properties listed in  \cite[Theorem 2.19]{\recognitionTag}.     Property (Z) applies only to zero strata and so is vacuous in this context.  Properties (F) and (NEG) are immediate from (a) - (c).    The endpoints of an \iNp\ are not contained in the interior of a fixed edge and so are vertices by Fact~\ref{ct basics}(1). This verifies the (V) property of \cite[Theorem 2.19]{\recognitionTag}.   If  a stratum $H_m$ is a forest in $\Per(f)$ then it is a single fixed edge $E_m$ with endpoints in a core subgraph $G_{m-1}$  by (a) - (c).  It follows that the free factor support of $G_{m-1}$ is not equal to the free factor support of $G_l \cup E_m$ for any filtration element $G_l$.  This verifies  the (P) property of \cite[Theorem 2.19]{\recognitionTag} and we are done.    
\endproof

 We assume for the rest of this subsection that 
 \begin{itemize}
 \item $\theta \in \upg$ and that $\fG$ is a \ct\ representing $\theta$, hence $f$ satisfies the conclusion of Lemma~\ref{ct basics}.   
\end{itemize}

If a root-free $a \in F_n$ is fixed by $m  \ge 2$ elements of $\P(\theta)$ then its unoriented conjugacy class $[a]_u$ is called an {\em axis} or {\em twistor} for $\theta$ with {\em multiplicity} $m-1$.   
 An edge $E$ in a stratum $H_i$  is {\em linear} if there is a Nielsen path $u \subset G_{i-1}$ such that $f(E) = Eu$.   Recall from Fact~\ref{boundary identification} that each non-trivial $a \in F_n$ corresponds to a covering translation $T_a : \ti G \to \ti G$ with axis $A_a$.

\begin{fact} \label{fact:twistor}  For each root-free $a \in F_n$,  if  $[a]_u$ is a twistor for $\theta$ of multiplicity $m \ge 2$ then :
\begin{enumerate}
\item There is a closed path $w$ that determines a  circuit representing $[a]_u$.
\item  There are exactly $m -1$ linear edges $E^1, \ldots, E^{m-1}$ such that $f(E^l) = E^lw^{d_l}$ for some  $d_l \ne 0$.  Furthermore, the  values  $d_l$ are pairwise distinct.  We say that $w$ is the {\em twist path} for  $E^l$.
\item For each $l=1,\ldots,m-1$, there is a lift $\ti f_l$ of $f$ such that for each lift $\ti E^l$, if the terminal endpoint of $\ti E^l$ is contained in $A_a$ then the initial endpoint of $\ti E^l$ is fixed by $\ti f_l$.  These lifts are pairwise distinct, preserve $A_a$  and each acts without fixed points on $A_a$.
\item   For each $l=1,\ldots,m-1$, the lift $\ti f_l$ corresponds to an element $\Theta_l \in \P(\theta)$ that fixes $a$.  These automorphisms $\Theta_l$ account for all but one element $\Theta_0 \in \P(\theta)$ that fixes $a$.  The lift $\ti f_0$ of $f$ that corresponds to $\Theta_0$ fixes points in $A_a$.
\end{enumerate}
\end{fact}

\proof  The (NEG Nielsen Paths) property of a \ct\ implies that  $\Fix(\ti f_i) \cap A_a = \emptyset$. The rest of the fact  follows from  the (Linear Edges)  
property of a \ct\ and  \cite[Lemma 4.40]{\recognitionTag}.
\endproof

\begin{fact}  \label{neutral crossing} If $E$ is a linear edge with twist path $w$ then every occurence of $E$ in a Nielsen path $\rho$ for $f$ is contained in a subpath of $\rho$ of the  form $E w^p \bar E$.   
 \end{fact}

\proof  We may assume that $\rho$ is indivisible.  Let $r$ be the height of $\rho$ and $s$ the height of $E$.  The $r < s$ case is vacuous and the $r=s$ case follows from the  (NEG Nielsen Paths) property for a \ct.  If $r> s$ then the  (NEG Nielsen Paths) property implies that $\sigma = E' {w'}^p \bar E'$  for some $p \ne 0$, where $E'$ is a linear edge of height $r$ and its twist path $w'$ has height less than $r$.  Each occurence of $E$ in $\sigma$ is contained in ${w'}^p$ and we are done by induction. 
\endproof

\begin{notation}\label{RsubE}    Letting $E$ be a non-fixed  edge of height $r$ with $f(E) = Eu$ for some closed path $u$ of height $< r$ [Fact~\ref{ct basics}\pref{item:edge in ct}],  its iterates split as $f^k(E) = E\cdot u \cdot f_\#(u)\cdot \ldots \cdot f_\#^{k-1}(u)$.   In this case, the nested sequence $E \subset  f(E) \subset f^2(E) \subset \ldots$ converges to {\em a ray $R_E$ that we say is  determined by $E$}.  If $E$ is a linear edge with twist path $w$ then $R_E = E w^{\pm \infty}$.   If $E$ is non-linear then the set of terminal endpoints of lifts of $R_E$ to $\ti G$ is  an $F_n$-orbit in $\partial F_n$ that we denote $[\partial R_E]$.  
 \end{notation}
    
\begin{fact}  \label{fact:eigenray} The assignment $E \mapsto [\partial R_E]$ defines  a bijection between the set  $\E$ of non-linear, non-fixed edges of $G$ and the set of eigenrays of  $\theta$ (Notation~\ref{notn:eigenray}).    
\end{fact}

\proof   This is contained in \cite[Fact 1.49]{\SubgroupsOneTag}; see also \cite[Lemma 4.36]{\recognitionTag}.
\endproof

\begin{corollary} \label{axes and eigenrays do not fill} If $\F \sqsubset \{[F_n]\}$ is a $\theta$-invariant  one-edge extension then $\F$ carries every twistor and eigenray for $\theta$. 
\end{corollary}

\proof  By \cite[Theorem 4.5]{\recognitionTag}, there exists a \ct\ $\fG$ representing $\theta$ in which $\F$ is represented by a filtration element $G_s$.  Each non-fixed edge $E$ above $G_s$ satisfies $f(E) = E \cdot u$ for some non-trivial path $u \subset G_s$.  The corollary therefore follows from Fact~\ref{fact:twistor} and Fact~\ref{fact:eigenray}.
\endproof

\begin{fact} \label{FixN not empty}   $\Fix_N(\Theta) \ne\emptyset$ for all $\Theta\in \Aut(F_n)$ representing $\theta$.  
\end{fact}

\proof  Let 
 $\ti f : \ti G \to \ti G$ be the lift of $f$ corresponding to $\Theta$.  If $\Theta$ is principal then $\Fix_N(\Theta) \ne \emptyset$ by definition.  We may therefore assume that $\Theta$ is not principal and hence by  Fact~\ref{ct basics}\pref{principal lifts} that $\ti f$  is fixed point-free.  In this case there is a path $\ti \sigma \subset \ti G$ such that $\ti \sigma \cdot \ti f_\#(\ti \sigma) \cdot \ti f_\#^2(\ti \sigma) \cdot \ldots$ converges to a point in $\Fix_N(\Theta)$.  The construction of $\ti \sigma$ is carried out in the proof of   \cite[Proposition 5.4.3]{\BookOneTag}.  A more directly quotable reference is \cite[Lemma 6.4]{FeighnHandel:ctAlgorithm}.
\endproof

\begin{fact} \cite[Proposition 4.44]{\BookTwoTag}  \label{upg restricts} If the conjugacy class of the free factor  $F$ is $\theta$-invariant then $\theta|F$ is \upg. 
\end{fact}

The following lemma is not known for elements of $\Out(F_n)$ that are not \upg.

\begin{lemma}   \label{principal symmetry}  $\Theta \in \P(\theta) \Longleftrightarrow \Theta^{-1} \in \P(\theta^{-1})$.
\end{lemma}

\proof  By symmetry, it suffices to assume that  $\Theta \in \P(\theta)$ and prove that  $\Theta^{-1} \in \P(\theta^{-1})$.  If the rank of $\Fix(\Theta)$ is at least two then this follows from the Definition~\ref{def:principal} and the fact that $\Fix(\Theta) = \Fix(\Theta^{-1})$.  We may therefore assume that $\Fix(\Theta)$ has rank one or zero. 
 We show below that   there is an injective map  $\Fix_+(\wh \Theta) \to  \Fix_+(\wh{\Theta^{-1}})$.     Assuming this for now, we complete the proof as follows. The cardinality of  $\Fix_+(\wh{\Theta^{-1}})$ is at least one in the rank one case and  at least two in the rank zero case.     It follows that $\Fix_N(\wh{\Theta^{-1}})$ contains  at least three points unless $\Fix(\Theta^{-1})$ has rank zero and $\Fix_+(\wh{\Theta^{-1}})$ contains exactly two points.  In this case, Fact~\ref{fix both} implies that $\Fix_+(\wh{\Theta^{-1}}) \ne \{a^\pm\}$ for any non-trivial $a \in F_n$ so $\Theta^{-1} \in \P(\theta^{-1})$ in this case as well.

 It remains to show that there is an injective map  $\Fix_+(\wh \Theta) \to  \Fix_+(\wh{\Theta^{-1}})$.  The lift  $\ti f : \ti G \to \ti G$  corresponding to $\Theta$ satisifies   $ \hat f = \wh \Theta$.  For each $P \in \Fix_+(\hat f)$ there is (Fact~\ref{fact:eigenray}) a non-fixed non-linear edge $E$, a lift $\ti E$ of $E$  and a lift $\wt{ R_E}$ of $R_E$ (Notation~\ref{RsubE}) with initial edge $\ti E$, such that $\wt{ R_E}$ converges to $P$  and  intersects $\Fix(\ti f)$ only in its initial endpoint $\ti v$ \cite[Lemma 3.36]{\recognitionTag}.  Let $r$ be the height of $E$.  By Lemma~\ref{ct basics}\pref{item:edge in ct}  there is a   component $C$ of $G_{r-1}$ that contains the terminal endpoint of $E$ and hence contains all of $R_E$ but its first edge.  Let $\Gamma \subset \ti G$ be the component of the full pre-image    $C$ that contains the terminal endpoint of the initial edge $\ti E$ of $\wt R_E$.    Then $\Gamma$ is $\ti f$-invariant and  $\partial \Gamma$ contains $P$. The (NEG Nielsen Paths) property of a \ct\ implies that $\ti f \restrict \Gamma$ is fixed point free. \cite[Lemma 3.16]{\recognitionTag} therefore implies that   $P$ is the only element of $\Fix_+(\hat f)$  contained in $\partial \Gamma$.  Since $\Fix_+(\hat f)$ is $\hat i_a$-invariant for each $a \in \Fix(\Theta)$, it follows that   
$\partial \Fix(\Theta) \cap \partial \Gamma = \emptyset$.

 Let $F$ be the  free factor that represents the unique element of $\F(C)$ and satisfies $\partial F  = \partial \Gamma$.      
The automorphism  $ \Psi:= \Theta^{-1} \restrict F$ represents the restriction  $\psi =\theta^{-1} \restrict F$, which  is \upg\ by Fact~\ref{upg restricts}.  By Fact~\ref{FixN not empty}   there exists  at least one point $Q \in \Fix_N(\Psi)$.  Since $\partial \Fix(\Psi) = \partial \Fix(\Theta) \cap \partial \Gamma= \emptyset$,  $Q \in \Fix_+(\wh \Psi)$. 

To see that $P \mapsto Q$ is injective, suppose that $P' \ne P$ is a point in $\Fix_+(\wh \Theta)$ and that $C', \Gamma', F'$ and $Q'$ are defined as above with $P$ replaced by $P'$.    Since $C$ and $C'$ are components of filtration elements of $G$, either they are disjoint or one is contained in the other.  It follows \cite[Fact 1.2]{\SubgroupsOneTag} that either $\partial F$ and $\partial F'$ are disjoint or one is contained in the other.  The latter is ruled out by the fact that   $P = \Fix_+(\Theta) \cap \partial F$ and $P' = \Fix_+(\Theta) \cap \partial F'$.   Thus  $\partial F$ and $\partial F'$  are disjoint and $Q \ne Q'$.  
\endproof

\begin{definition} \label{def:highest edge splitting}
Every path $\ti \sigma \subset \ti G$ with endpoints, if any, at vertices has a {\em highest edge splitting} $\ti \sigma = \ldots \ti \sigma_{-1} \cdot \ti \sigma_0 \cdot \ti \sigma_1 \ldots$ defined as follows.   If $r$ is the height of $\sigma$ and $E_r$ is not fixed  then this  splitting is  defined by taking the splitting vertices ( i.e. the endpoints of the terms)  to be exactly those vertices that are either the initial endpoint of an edge in $\ti \sigma$ that projects to $E_r$ or the terminal endpoint of an edge that projects to $\bar E_r$.  If $E_r$ is fixed then both endpoints of an edge that projects to $E_r$ or $\bar E_r$ are splitting vertices. The projected splitting $\sigma = \ldots  \sigma_{-1} \cdot \sigma_0 \cdot  \sigma_1 \ldots$ is the {\em highest edge splitting} of $\sigma$.
\end{definition}

\begin{fact} \label{fact:highest edge splitting}
 \begin{enumerate} 
 \item The highest edge splitting $\sigma = \ldots  \sigma_{-1} \cdot \sigma_0 \cdot  \sigma_1 \ldots$ of $\sigma$ is in fact a splitting. 
 \item \label{item:hes} For any lift $\ti f$, $$\ti f_\#(\ti \sigma)  = \ldots \ti f_\#(\ti \sigma_{-1}) \cdot \ti f_\#(\ti \sigma_0) \cdot \ti f_\#(\ti \sigma_1) \ldots$$ is the highest edge splitting of $\ti f_\#(\ti \sigma)$. 
 \end{enumerate}
\end{fact}

\proof  Item (1) is contained in the statement and proof of \cite[Lemma~4.1.4]{\BookOneTag}.  

   For \pref{item:hes}, let $\V_{\ti \sigma}$ and $\V_{\ti f_\#(\ti \sigma)}$  be the  highest edge splitting   vertices of $\ti \sigma$ and $\ti f_\#(\ti\sigma)$ respectively.   Assuming at first that $E_r$ is not fixed,   each term $\ti \sigma_j$  in the  highest edge splitting of $\ti\sigma$ have the form  $\ti E_r\ti \gamma  \ti E_r^{-1}$, $\ti E_r \ti \gamma$, $\ti \gamma   \ti E_r^{-1}$, $\ti \gamma$, $\ti E_r$ or $\ti E_r^{-1}$ for some non-trivial path $\ti\gamma$ that projects into $ G_{r-1}$.   Since  $f(E_r) =  E_r u_r$ for some path $u_r \subset G_{r-1}$, the $\ti f_\#$-image of each of these types is another path of the same type. It follows that  $\V_{\ti f_\#(\ti \sigma)} \subset \ti f(\V_{\ti \sigma})$.    It also follows that if $\ti \sigma_j$ ends with $\ti E_r^{-1}$ [respectively begins with $\ti E_r]$ then $\ti f_\#(\ti \sigma_j)$ ends with $\ti E_r^{-1}$  [respectively begins with $\ti E_r]$.   This implies that $ \ti f(\V_{\ti \sigma}) \subset \V_{\ti f_\#(\ti \sigma)} $.  This completes the proof in the case that $E_r$ is not fixed.  The remaining case is similar and is left to the reader. 
\endproof

\begin{lemma} \label{line without fixed points}  Suppose that $\ti f$ is a lift of $\fG$, that $\ti \mu$ is an $\ti f_\#$-invariant line that is disjoint from $\Fix(\ti f)$ and that an endpoint of $\ti \mu$ is fixed by a covering translation $T$.  Then $\ti \mu$ is the axis $A_T$ of $T$.
\end{lemma}

\proof    Let $\ti \mu = \ldots \cdot\ti \mu_{-1} \cdot \ti  \mu_0 \cdot\ti  \mu \cdot \ldots$ be  the highest edge splitting of $\mu$. Fact~\ref{fact:highest edge splitting} implies that  there exists $p \in \Z$ such that $(\ti f)_\#(\ti  \mu_i) = \ti  \mu_{i+p}$ for all $\ti  \mu_i$.   From the assumption that $\Fix(\ti f) \cap \ti \mu = \emptyset$, it follows that $p \ne 0$ and so the splitting is  bi-infinite.   The highest edge splitting $A_T =  \ldots \cdot \ti \alpha_{-1}  \cdot\ti \alpha_0 \cdot \ti \alpha_1  \cdot \ldots$  of $A_T$ is also bi-infinite.  Since   $\ti \mu$ and $A_T$  have a common ray, they must have the same height.  After re-indexing the $\ti \alpha_j$'s, we may assume that  $ \ti \alpha_j = \ti  \mu_j$ for all sufficiently large $j$.  It follows that   $(\ti f)_\#(\ti \alpha_j) = \ti \alpha_{j+p}$ for all $\ti \alpha_j$  and hence   that for all $j$ there exists $k > 0$ such that  $(\ti f)^k_\#(\ti  \mu_j)  = (\ti f)^k_\#(\ti \alpha_j)$.  Since $ \mu_j$ and $\alpha_j$ are paths in $G$ with the same endpoints, it follows that $ \ti \alpha_j = \ti  \mu_j$ for all $j$ and  $\ti \sigma = A_T$. 
\endproof

 \subsection{Abelian \upg\ subgroups}
  {\em We assume throughout this section that $K$ is an abelian \upg\ subgroup  of $\Out(F_n)$.}   Lemma~\ref{upg is rotationless} implies that each element of $K$ is rotationless and so $K$ is a rotationless abelian \upg\ subgroup.    
    
The main definitions in \abelian\ make use of    
$$\P^\pm(\theta) := \P(\theta) \cup (\P(\theta^{-1}))^{-1}$$
 In the \upg\ case,   Lemma~\ref{principal symmetry} implies that 
 $$\P^\pm(\theta) = \P(\theta) = \P(\theta^{-1})^{-1}$$
We have simplified the definitions in this subsection accordingly.      
    
\begin{remark} As noted in \cite[Section 6.1.2]{HandelMosher:BddCohomologyII}, the definition of $P^{\pm}(\theta)$ was misstated in \recognition\ as $\P^\pm(\theta) := \P(\theta) \cup \P(\theta^{-1})$.
\end{remark}  
   
 \begin{definition} \cite[Definition 3.9]{\abelianTag} \label{principal lift} \ \ A set $X \subset \partial F_n$ with at least three points is a {\em principal set} for $K$ if for each $\phi \in K$ there exists $\Phi \in \P(\phi)$ such that $X \subset \Fix(\wh \Phi)$.  For each such $X$, the assignment $\phi \to \Phi$ defines a lift $s_X$ of $K$ into $\Aut(F_n)$ called the {\em principal lift} determined by  $X$.  
 \end{definition}
 
 If $X$ is a principal set  and $X' \subset X$  contains at least three points then $X'$ is a principal set and    $s_{X'} = s_X$.   For any principal set $X'$, the maximal (with respect to inclusion) principal set  containing $X'$ is given by $$X = \cap_{\phi \in K} \Fix(\wh{s_{X'}(\phi})$$ See \cite[Remark 3.10]{\abelianTag}.   
 
   Automorphisms that differ by conjugation by an inner automorphism are said to be {\em isogredient}.  If $X$ is a maximal principal set and $c \in F_n$ then $i_c(X)$ is also a maximal principal set and $s_{i_c(X)}(\phi) = i_c s_X (\phi)  i_c^{-1}$ for each $\psi \in K$. Thus $F_n$-orbits of  maximal principal sets correspond to isogredience classes of principal lifts.   
   
   The following definition generalizes to the setting of \upg\ abelian subgroups the concepts that were defined just preceding Fact~\ref{fact:twistor}.
 
\begin{definition}   \cite[Definition 4.1]{\abelianTag} \label{A(K)}  We say that the unoriented conjugacy class $[a]_u$ of $a \in F_n$ is an {\em axis}  or {\em twistor} of multiplicity $m-1\ge1$  for $K$ and write $[a]_u \in \A(K)$  if $\{a^\pm\}$ is  contained in $m$ distinct maximal principal sets.   The maximal  principal sets that contain $\{a^\pm\}$ are called  {\em linear principal sets} or more specifically  {\em $a$-linear principal sets}.     If $X_1$ and $X_2$ are distinct $a$-linear principal sets then for each $\theta \in K$ there exists an integer $d(\theta)$ such that   $s_{X_2}(\phi) = i_a^{d(\theta)} s_{X_1(\theta)}$.  The assignment $\theta \mapsto d(\theta)$ defines a homomorphism $\omega : K \to \Z$ called the {\em comparison homomorphism  determined by $X_1$ and $X_2$}.    Note that $\omega$ depends only on the $F_n$-orbit of the pair $(X_1,X_2)$; i.e. $(X_1,X_2)$ and $i_c(X_1,X_2) := (i_cX_1,i_cX_2)$ determine the same comparison homomorphism.   
\end{definition}

\begin{fact} \cite[Lemma 4.3]{\abelianTag}  There are only finitely many comparison homomorphisms for $K$.
\end{fact}
\begin{fact} \label{coordinates} \cite[Lemma 4.6]{\abelianTag}  If $\theta, \phi \in K$ and $\omega(\theta)= \omega(\psi)$ for all comparison homomorphisms $\omega$ then $\theta = \phi$.
\end{fact}

\begin{definition} \cite[Definition 4.7]{\abelianTag}  {\em $\phi \in K$ is generic} if $\omega(\phi) \ne 0$ for each comparison homomorphism $\omega$.
\end{definition}

 \begin{fact}\cite[Lemma 4.10]{\abelianTag}  \label{generic fixed sets} If $\theta \in K$ is generic then $\{\Fix(\wh{\Theta}): \Theta \in \P(\theta)\}$ is the set of maximal principal sets for  $K$. 
\end{fact}

\begin{cor}  \label{same axes}If $\phi, \theta \in K$ are generic then $\{\Fix(\wh{\Theta}): \Theta \in \P(\theta)\} = \{\Fix(\wh{\Phi}): \Phi \in \P(\phi)\}$.
\end{cor}

\begin{cor}  \label{generic twistor} $[a]_u \in F_n$ is  a twistor for $K$ of multiplicity $m-1$  if and only if $[a]_u$ is a twistor of multiplicity $m-1$ for some, and hence every, generic element of $K$.
\end{cor}

\begin{fact}\label{basis of generics} \cite[Lemma 4.9]{\abelianTag}  $K$ has a basis of generic elements.
 \end{fact}

\begin{fact}\label{conjugation} \label{conjugating P}\cite[Lemma 2.6]{\abelianTag}   Suppose that  $\theta, \psi  \in  \Out(F_n)$, that $\phi := \theta^\psi= \psi \theta \psi^{-1}$ and that $\Psi \in \Aut(F_n)$ represents   $\psi $. 
 Then 
\begin{enumerate}
\item $\Fix(\wh{\Psi \Theta \Psi^{-1}}) = \wh \Psi(\Fix\wh \Theta)$  for all $\Theta \in \Aut(F_n)$ representing  $\theta$.
\item  $\Fix_+(\wh{\Psi \Theta \Psi^{-1}})=\wh \Psi(\Fix_+(\wh \Theta))$ and $ \Fix_-(\wh{\Psi \Theta \Psi^{-1}}) = \wh \Psi(\Fix_-(\wh \Theta)) $ for all $\Theta \in \Aut(F_n)$ representing  $\theta$.
.
\item $\Theta \mapsto  \Psi \Theta \Psi^{-1}$ defines a bijection between $\P(\theta)$ and $\P(\phi)$ that preserves isogredience classes.

\end{enumerate}
\end{fact}

\begin{lemma} \label{still generic}   If $\theta$ is generic in  $K$  and $\psi \in \Out(F_n)$ then $\theta^\psi  $ is generic in $K^\psi = \{\psi \phi \psi^{-1}: \phi \in K\}$.  Moreover, for any $\Psi \in \Aut(F_n)$ representing $\psi$, $\wh \Psi$ induces a bijection between [$a$-linear] principal sets in $K$ and [$\Psi(a)$-linear] principal sets in $\K^\psi$.
\end{lemma}

\proof   Choose $\Psi \in \Aut(F_n)$ representing $\psi$.  The following are easy consequences of Fact~\ref{conjugation}:
\begin{itemize}
\item  $\A(K^\psi) = \psi(\A(K))$.
\item If $X_i$ is a  [$a$-linear] principal set  for $K$ then $\wh\Psi(X_i)$ is a [$\Psi(a)$-linear] principal set for $K^\psi$. 
\item If $\omega : K \to \Z$ is the comparison homomorphism determined by $X_1$ and $X_2$ then the  comparison homomorphism $\omega^\psi : K^\psi \to \Z$ determined by $\wh \Psi(X_1)$ and $\wh \Psi(X_2) $ satisfies $\omega^\psi(\theta^\psi) =    \omega(\theta)$.  
\end{itemize}
The lemma now follows from the definition of genericity.
\endproof

\begin{lemma}  \label{generic restriction}   If   $\theta \in K$ is generic in $K$ and $[F]$ is a $K$-invariant free factor conjugacy class  then $K \restrict F$ (see Fact~\ref{restriction}) is an  abelian \upg\ subgroup and 
$\theta \restrict F$ is generic in $K \restrict F$. 
\end{lemma} 

\proof  Fact~\ref{upg restricts} implies that $K \restrict F$ is an  abelian \upg\ subgroup. Each maximal principal set for $K \restrict F$ extends uniquely to a maximal principal set for $K$.  It follows that each coordinate homomorphism $\omega_K$  for $K \restrict F$ is the restriction of a coordinate homomorphism $\omega$ for $K$. Thus each $\omega_K(\theta \restrict K) = \omega(\theta) \ne 0$ proving that $\theta \restrict K$ is generic.   
   \endproof

\section{Proof of Proposition~\ref{upg prop}} \label{main proof} 
Assuming that  $K \subgroup \IAThree$ is an  abelian \upg\ subgroup and  that $\psi \in \IAThree$ normalizes $K$, our goal is to  show that $\psi$ commutes with each $\theta \in K$. By Facts~\ref{upg is rotationless} and  \ref{basis of generics},   $K$ is rotationless  and has a basis of   of generic elements.  We may therefore   assume that $\theta$ is  generic.   Letting $\phi = \theta^\psi =  \psi \theta \psi^{-1} \in K$,  
 our goal is to show that $\phi = \theta$.

\medskip

 Suppose that $[a]_u \in \A(K)$ (see Definition~\ref{A(K)}) and that $\Psi$ is a representative of $\psi$.  Lemma~\ref{still generic} implies that $[\Psi(a)]_u \in \A(K^\psi) = \A(K)$ and hence that $\psi$ permutes the elements of $\A(K)$.  It follows that $\psi$ fixes each element of $\A(K)$ by   \cite[Theorem 4.1]{\SubgroupsTwoTag}.   We can therefore choose $\Psi_a$ representing $\psi$  that fixes $a$.       Lemma~\ref{still generic} implies that   if $X_i$ is an $a$-linear principal set for $K$ then 
$\wh{\Psi_a}(X_i)$ is an $a$-linear principal set for $  K^\psi = K$.    Thus
\begin{description}
\item [($*$)] For each $a \in F_n$ with $[a]_u \in \A(K)$ there exists $\Psi_a \in\Aut(F_n)$ representing $\psi$ such that $\Psi_a(a) =a$. Letting $\{X_0,\ldots, X_{m-1}\}$ be the  $a$-linear principal sets, the automorphism $\Psi_a$     induces a permutation $\pi_a$ of $ \{0,\dots,m-1\}$ such that  $\wh{\Psi_a}(X_i) = X_{\pi_a(i)}$.   Corollary~\ref{generic twistor} implies  that $[a]_u$ is   twistor for $\theta$ with multiplicity $m-1$.  
\end{description} 
The main work of the proof is to  show  that each $\pi_a$ is the identity.   Assuming this fact for the moment,  we complete the proof of the proposition.

\medskip
Fix $[a]_u \in \A(K)$ and let $\{X_0,\ldots, X_{m-1}\}$, $\Psi_a$ and $\pi_a$ be as in ($*$).
Let $\Theta_i = s_{X_i}(\theta) \in \P(\theta)$ (see Definition~\ref{principal lift}) and $\Phi_i =  \Psi_a \Theta_i \Psi_a^{-1} \in \P(\phi)$.  Then 
$$\Fix(\wh \Phi_i) = \wh{\Psi_a}(\Fix(\wh \Theta_i)) = \wh {\Psi_a}(X_i) = X_{\pi_a(i)} = X_i$$
 with the middle equality following from Fact~\ref{generic fixed sets} and the genericity of $\theta$.  Thus $s_{X_i}(\phi) = \Phi_i$.      If $X_i$ and $X_j$ are distinct $a$-linear principal sets then  
$$\Theta_j =  i_a^d \Theta_i$$ for some $d \ne 0$ and
$$\Phi_j  =\Psi_a \Theta_j \Psi_a^{-1} = \Psi_a i_a^d \Theta_i  \Psi_a^{-1} = i_{\Psi_a(a)}^d \Psi_a  \Theta_i  \Psi_a^{-1} = i_a^d \Phi_i$$
This proves that $\omega(\theta) = \omega(\phi)$ where $\omega$ is the comparison homomorphism  determined by $X_i$ and $X_j$.  Since $a,i$ and $j$  are arbitrary, Fact~\ref{coordinates}  completes the proof of Proposition~\ref{upg prop}.  

\medskip

Fixing $a$ as in ($*$), it remains to prove that $\pi_a$ is the identity.  

\medskip

\begin{lemma}  \label{axes and eigenrays} Suppose that $\emptyset = \F_0 \sqsubset\F_1 \sqsubset \ldots \sqsubset \F_{J} = \{[F_n]\}$  is a maximal nested sequence of free factor systems  that are invariant by both   $K$ and $\psi$ and  that     $[F]$ is a component of $\F_j$ for some $1\le j \le J$.  Then $[F]$  is invariant by both   $K$ and $\psi$; moreover,    for all $\eta \in K$,   each axis and eigenray  for  $\eta \restrict F$   is carried by $\F_{j-1}$.
\end{lemma}

\proof  Since $K, \psi \subset \IAThree$,  \cite[Lemma 4.2]{\SubgroupsTwoTag} implies that $[F]$ is invariant by both   $K$ and $\psi$.    By Fact~\ref {upg restricts},  $K_F:= K\restrict F$ is a  \upg\ subgroup   that is obviously normalized by  $\psi_F := \psi \restrict F$.  Each component of $\F_{j-1}$ is contained in a unique component of $\F_j$.  The union of the components of $\F_{j-1}$ that are contained in $[F]$ define a free factor system of $F$ that we denote by   $\F'$.   Since $\F_{j-1} \sqsubset \F_{j}$ is   invariant by both   $K$ and $\psi$ and is  maximal   with respect to these properties, it follows that $\F' \sqsubset \{[F]\}$ is invariant by both   $K_F$ and $\psi_F$ and is  maximal   with respect to these properties.    By  \cite[Theorem 5.1]{\BookTwoTag} there is a $K_F$-invariant free factor system $\F''$ of $F$ such that $\F'\sqsubset \F''$ and such that $\F''\sqsubset \{[F]\}$ is a one-edge extension.  By Corollary~\ref{axes and eigenrays do not fill},  $\F''$ carries  each axis and eigenray  for all $\eta_F \in K_F$.    It follows from Fact~\ref{conjugating P} and the definitions that $\psi$ maps the axes and eigenrays of $\eta_F \in K_F$ to the axes and eigenrays of $\eta_F^{\psi_F} \in K_F$ and so the set $X$ of axes and eigenrays that occur for some element of $K_F$ is $\psi_F$-invariant.  For the same reason, $X$ is $K$-invariant.  By Fact~\ref{ffs support}, $\F_\supp(X \cup \F')$  is invariant by both $K$ and $\psi$. Moreover, $\F_\supp(X \cup \F')$ contains $\F'$ by construction and   is properly contained in $\{[F]\}$  because it is contained in $\F''$.  It follows that  $\F_\supp(X \cup \F') = \F'$ and so  $\F_\supp(X ) \sqsubset \F' \sqsubset \F_{j-1}$.
\endproof
  
 Having fixed $a$ as in  ($*$), the notations of ($*$) remain in force for   the rest of the paper, as do Notations ~\ref{standing notation} and \ref{Sigma} below.  Recall also that $\theta$ is a generic element of $K$, that $\psi$ normalizes $K$ and that $\phi= \theta^\psi$.
\begin{notation}\label{standing notation}Let $\emptyset = \F_0 \sqsubset\F_1 \sqsubset \ldots \sqsubset \F_{J} = \{[F_n]\}$  be a maximal nested sequence of free factor systems  that are invariant by both   $K$ and $\psi$.  By  \cite[Theorem 1.1]{FeighnHandel:ctAlgorithm} there is a \ct\    $\fG$ representing $\theta$ with filtration $\emptyset = G_0 \subset G_1 \subset \ldots  \subset G_N = G$    and a subfiltration by core graphs $\emptyset = G_0 = G_{r(0)} \subset G_{r(1)} \subset \ldots\subset G_{r(J)}=G$  such
each $\F_j$ is realized by  $G_{r(j)}$; moreover, for each $j$ and each component  $C$ of $G_{r(j)}$, \  $f|C$ is a \ct\ representing $\theta \restrict [\pi_1(C)]$.  For each $1 \le j \le J$,  define the $j^{th}$-stratum of  the {\em subfiltration} to be the subgraph $S_j = G_{r(j )}\setminus G_{r(j-1)}$.  Choose $h : G\to G$ representing $\psi$ such that each $G_{r(j)}$ is $h$-invariant.   
\end{notation}

 \begin{cor} \label{high enough edges} If $E_i$ is a non-fixed edge in $S_j$  then $f(E_i) = E_i \cdot u_i$ for some non-trivial closed path $u_i \subset G_{r(j-1)}$.
\end{cor}

\proof  Let $C$ be the component of $G_{r(j)}$ that contains $E_i$ and let $F$ be a free factor representing $[\pi_1(C)]$.  
 By construction, $f \restrict C$ is a \ct\ representing $\theta \restrict F$. If $E_i$ is linear then $u_i = w_i^{p_i}$ where $w_i$ is the twist path for $E_i$ and $[w_i]_u$  is an axis for $\theta \restrict F$.  In this case,  Lemma~\ref{axes and eigenrays} completes the proof.   If $E_i$ is non-linear then $\partial R_{E_i}$ is an eigenray for $\theta \restrict F$ by Fact~\ref{fact:eigenray}.  Lemma~\ref{axes and eigenrays} therefore implies that $f^k_\#(u_i)  \subset G_{r(j-1)}$ for some $k > 0$ which implies that $u_i  \subset G_{r(j-1)}$.  
\endproof

\begin{notation} \label{Sigma}For $0 \le i \le m-1$, let $\Theta_i = s_{X_i}(\theta)$,  let $\ti f_i : \ti G \to \ti G$ be the lift corresponding to $\Theta_i$ and let $\hat f_i:\partial F_n \to \partial F_n$ be the extension of $\ti f_i$ given by Fact~\ref{f hat}.  Then $\Theta_0,\ldots, \Theta_{m-1}$ are  the elements of $\P(\theta)$ that fix $a$  and $\Fix(\hat f_i) = X_i$ by Facts~\ref{boundary identification} and \ref{generic fixed sets}.
 After possibly reindexing the $X_i$'s,  we may assume that $\Theta_1,\ldots, \Theta_{m-1}$ correspond to linear edges $E^1,\ldots, E^{m-1}$ as in Fact~\ref{fact:twistor}.       We may also assume   that the twist path for $[a]_u$    represents $[a]$.  Superscript indices, primarily  $i,j,k$ or $l$, take values in the set $\{1,\ldots,m-1\}$.   Edges indexed by subscripts like $E_p$ or $E_q$ can be any edge at all, perhaps even an element of $\{E^1,\ldots, E^{m-1}\}$. 
 
  The twist path for  $E^1,\ldots, E^{m-1}$ is denoted by $w$.   We define rays by iterating $w$ and $w^{-1}$ in the positive and negative directions as follows:
    $$R_+(w) = (www\ldots) \qquad 
    R_+(\bar w) = (\bar w \bar w \bar w\ldots) $$
    $$ 
    R_-(w) = (\ldots www) \qquad
    R_-(\bar w) =(\ldots \bar w \bar w \bar w)$$
Given a line $\alpha$, we say that {\em $\alpha$   ends with $w^\infty$} if there is a concatenation expression $\alpha = \beta R_+(w)$. Similarly, {\em $\alpha$  begins with $w^\infty$}  if $\alpha = R_-(w) \beta$; $\alpha$ {\em ends with $\bar w^\infty$} if $\alpha = \beta R_+(\bar w)$; and {\em $\alpha$ begins with $\bar w^\infty$} if $\alpha = R_-(\bar w) \beta$. Since the notation should make the context clear, we will usually abuse notation by ignoring \lq R\rq\ and writing $\beta w^\infty$ instead of $\beta R_+(w)$, and similarly for the other three possibilities.
  \end{notation}

\begin{definition} Recall that all lines in this paper are oriented. Since $X_i$ contains $a^+, a^-$ and at least one other point, there exist lines in $\ti G$  with initial endpoint in $X_i \setminus \{a^\pm\}$ and terminal endpoint in $\{a^\pm\}$.  For $0 \le i \le m-1$, let $\Sigma_i$  be the set of such lines of minimum \lq subfiltration height\rq\  $j(i) \in \{1,É,J\}$. To be more precise, let $1 \le j(i) \le J$ be the minimum value    for which there  exists an $(\ti f_i)_\#$-invariant  line $\ti \sigma \ne A_a$ (equivalently, a  line $\ti \sigma \ne A_a$ with endpoints in $X_i$) that terminates at either $a^+$ or $a^-$ and whose projection $\sigma$ is contained in $G_{r(j(i))}$.  The set of all such $\ti \sigma$ is denoted by   $\wt{\Sigma_i}$  and the set of  projections $\sigma$ is denoted by $\Sigma_i$.  Note that every $\sigma \in \Sigma_i$ ends with $w^\infty$ or $w^{-\infty}$.
\end{definition} 

 \begin{remark}  If $i \ne 0$ then  by Lemma~\ref{splitting at fixed points}\pref{item:last term infinite} below, $j(i)$ is the minimum value   for which $G_{r(j(i))}$ contains $E^i$.  
 \end{remark}
 
 \begin{remark}  \label{plus minus} If $\sigma \in \Sigma_i$  decomposes as $\sigma =\beta w^{-\infty}$ then there is a lift $\ti \sigma$ with terminal endpoint $a^-$ and initial endpoint, say $P \ne a^+$, in $ X_i$.  The line $\ti \sigma'$ that has initial endpoint $P$ and terminal endpoint $a^+$ projects to $\sigma' \in \Sigma_i$ that  decomposes as $\sigma' = \beta' w^{\infty}$.     
 \end{remark}

\begin{lemma} \label{Fsupp}If $\pi_a(k) = i$ then:
\begin{enumerate}
\item  $j(k)=j(i)$.
\item \label{item:action of psi}$\psi_\#(\Sigma_k) = \Sigma_i$.
\item \label{item:supp Sigma}$\F_\supp(\Sigma_k) = \psi(\F_\supp(\Sigma_k)) = \F_\supp(\Sigma_i)$. 
\end{enumerate}
\end{lemma}

\proof     
Recall that $h :G \to G$ represents $\psi$ and that each $G_{r(j)}$ is $h$-invariant.   Let $\ti h : \ti G \to \ti G$ be the lift corresponding to $\Psi_a$.  In particular, $\hat h = \hat \Psi_a$ fixes both $a^-$ and $a^+$.   If $\ti \sigma \in \Sigma_k$ then $\ti h_\#(\ti \sigma)$ has endpoints in $\wh \Psi_a(X_k) = X_i$ and   projects to $h_\#(\sigma) \subset G_{r(j(k))}$.  Thus $j(i) \le j(k)$.  This can be repeated  to show that $j(k) \ge j(\pi_a(k)) \ge j(\pi_a^2(k)) \ge \ldots$.  Since $\pi_a$ has finite order it follows that $j(\pi_a^l(k))$ is independent of $l$ and in particular that $j(k) = j(i)$.  Since $\hat \Psi_a (X_k \cap \partial G_{r(j(k))}) =   X_i \cap \partial G_{r(j(i))}$,  it follows that $h_\#(\Sigma_k) = \Sigma_i$ and hence that $\psi(\F_\supp(\Sigma_k)) = \F_\supp(\Sigma_i)$.  Iterating this argument, shows   that $\psi^l(\F_\supp(\Sigma_k)) = \F_\supp(\Sigma_{\pi_a^l(k)})$ for all $l \ge 1$.  Since $\pi_a$ has finite order,  $\F_\supp(\Sigma_k)$ is preserved by an iterate of $\psi$   and so  \cite[Lemma 4.2]{\SubgroupsTwoTag}  is  also preserved by $\psi$.   \endproof

Our strategy for proving that $\pi_a$ is the identity is to assume that this is not the case, and to produce a closed path $\delta \subset G$ representing a homology class in $H_1(G;\Z/3)$ that is not fixed by $\phi$. Typically $\delta$ will occur as a subpath of some nonperiodic line $\tau$ in some $\Sigma_i$ having the form $\tau = w^\infty \delta w^\infty $.   Homology information of such paths $\delta$ will be extracted from the algebraic crossing number of each line in each $ \Sigma_i$ with certain edges $E_p$. The information about lines in the sets $\Sigma_i$ that we need for these purposes is contained in the following lemma.Ó

\begin{lemma}  \label{splitting at fixed points} For $0 \le i \le m-1$ and $\ti \sigma \in \wt{\Sigma_i}$, let $\ti \sigma = \ldots \cdot \ti \sigma_{-1} \cdot \ti \sigma_0 \cdot\ti \sigma_1 \cdot \ldots$ be the decomposition where the endpoints of the $\ti \sigma_m$'s are exactly the   vertices of $\Fix(\ti f_i) \cap \ti \sigma$.
\begin{enumerate}  
\item \label{item:non-trivial splitting} The decomposition is a non-trivial splitting whose finite terms are fixed edges and indivisible Nielsen paths for $\ti f_i$.  
 \item  \label{item:finite term} If a finite $\ti \sigma_s$   crosses $E_p$ or $\bar E_p$ for some  non-fixed edge $E_p \subset S_{j(i)}$  then $E_p$ is linear and  $\sigma_s = E_p w_p^* \bar E_p$ where $w_p$ is the  twist path for $E_p$.   Moreover, $w_p^ {\pm \infty} \bar E_p\cdot \sigma_{s+1}\sigma_{s+2}\cdots$ is an element of $\Sigma_i$. 
\item \label{item:last term infinite} If $i \ne 0$ then the splitting has  a last term and it projects to $E^iw^{\pm \infty}$.
\item \label{item:no last term} If $i =0$ then the splitting has no last term.
\item \label{item: first term}  If the splitting has  a first term $\ti \sigma_0$    and if $\sigma_0$  crosses $E^t$ or $\bar E^t$ where  $E^t\subset S_{j(i)}$ and  $1 \le t \le m-1$ then   $\sigma_0$ has the form $w^{\pm \infty} \bar E^t$. 

\end{enumerate}
\end{lemma}

\proof   Non-triviality of the decomposition  follows from  Lemma~\ref{line without fixed points}.   It is a splitting because each of its terms is fixed by $\ti {f_i}_\#$.   
  For the same reason, each finite $\ti\sigma_l$ is a Nielsen path. If  $\ti\sigma_l$ is neither a single edge nor an \iNp\  then it would contain a fixed point in its interior;  the  (Vertices) property of a \ct\ would then imply that  $\ti\sigma_l$ contains a fixed vertex in its interior which is not the case.    This proves    \pref{item:non-trivial splitting}.

  If $\ti \sigma_s$ is finite but not a fixed edge then   
    $\sigma_s= E_q w_q^* \bar E_q$ for some linear edge $E_q \subset G_{r(j(i))}$ with twist path $w_q$ by \pref{item:non-trivial splitting} and the (\noneg\ Nielsen path) property of a \ct.   Corollary~\ref{high enough edges}  implies that $w_q \subset G_{r(j(i-1))}$ and  so does not cross any edge  in $S_{j(i)}$.  Thus $E_q = E_p$  and the main statement of  \pref{item:finite term}  is satisfied.  The moreover part of  \pref{item:finite term} follows from the following observations:  the turn $(E_p,\sigma_{s+1})$ is legal;  $\ti \sigma_{s+1} \cdot \ti\sigma_{s+2}\cdots$ and $ \ti  E_p \ti w_p^ {\pm \infty}$ are   $\ti {f_i}_\#$-invariant rays;   $\ti \sigma_{s+1} \cdot \ti\sigma_{s+2}\cdots$ and $ \ti  E_p \ti w_p^ {\pm \infty}$  project into $S_{j(i)}$.

For \pref{item:last term infinite} we assume that $i \ne 0$.    Lemma~\ref{fact:twistor}  implies that  $\Fix(\ti f_i) \cap A_a = \emptyset$ so there is a  last  fixed point in $\ti \sigma$ and a last term, say $\ti \sigma_b$. Since $\ti \sigma_b$ is an $\ti {f_i}_\#$-invariant ray, Fact~\ref{fact:highest edge splitting} implies that each term in the highest edge splitting of $\ti \sigma_b$ is  $\ti {f_i}_\#$-invariant.  If this splitting is non-trivial, the terminal endpoint of its first term   would be fixed in contradiction to the fact that the interior of $\ti \sigma_b$ is disjoint from $\Fix(\ti f_i)$.   The highest edge splitting  of $\ti \sigma_b$ is therefore trivial which means 
that   $\ti \sigma_b = \ti E_q \ti \rho$ where  $\ti E_q$ is   a  non-fixed edge with fixed initial direction and $\rho$ has height strictly less than that of  $E_q$.   Let $ R_{E_q}$ be the ray determined by $E_q$  (Notation~\ref{RsubE}).   The lift $\wt{R_{E_q}}$ that begins with $\ti{E_q}$ terminates at some 
$Q \in \partial \Fix(\Theta_i) \cup \Fix_+(\hat \Theta_i) \subset   X_i$ (see Fact~\ref{gjll}) and   intersects $\Fix(\ti f_i)$ only in its initial endpoint by (NEG Nielsen Paths).  
If $Q$ is not equal to the terminal endpoint $P$ of $\ti \sigma$ then the line $\ti L$ connecting $Q$ to $P$ is  disjoint from $\Fix(\ti f_i)$ and so equals $A_a$ by Lemma~\ref{line without fixed points}.  In this case,  $Q$ is either $a^+$ or $a^-$.  The same is true if  $Q = P$ because $P$ is either $a^+$ or $a^-$.  If $E_q$ is non-linear than $Q \in \Fix_+(\ti f_i)$ by 
Facts~\ref{fact:eigenray} and \ref{gjll} therefore imply that   $E_q$ is linear.  It follows that $\ti R_{E_q} \setminus \ti E_q$ is contained in   the axis of a covering translation that shares a terminal ray with $A_a$ and so equals $A_a$.  Since $\ti R_{E_q} \setminus \ti E_q$  projects to $w_q^{\pm \infty}$, $[w_q]_u = [w]_u$ and so $w_q = w$ and  the terminal endpoint of $\ti E_q$ is in $A_a$.  Combining this with Fact~\ref{fact:twistor} and the fact that the initial endpoint of $\ti E_q$ is fixed, we see that   that $E_q = E^i$.   This completes the proof of \pref{item:last term infinite}.

 If $i = 0$ then   $A_a \cap \Fix(\ti f_i)$ is non-empty by Lemma~\ref{fact:twistor} and is invariant under the covering translation $T_a$ associated to $a$   because $\ti f_i$ fixes $a^+$ and $a^-$.  It follows that $a^+$ and $a^-$ are in the closure of $A_a \cap \Fix(\ti f_i)$.    This proves \pref{item:no last term}. 
 
 The proof of \pref{item: first term} is similar to that of \pref{item:last term infinite}.   Assuming that the splitting has a first term $\tilde\sigma_0$ crossing $E^t$ or $\bar E^t \subset S_{j(i)}$, let $\ti R = \ti \sigma_0^{-1}$ and let $P$ be its terminal endpoint.  As in the proof of \pref{item:last term infinite}, the highest edge splitting of $\ti R$ must  be trivial   so  $ \ti R = \ti E_q \ti \rho$ where  $\ti E_q$ is   a  non-fixed edge with fixed initial direction and $\rho$ has height strictly less than that of  $E_q$.   
   Let $R_{E_q}$ be the ray determined by $E_q$, let $\wt R_{E_q}$ be the lift that begins with $\ti E_q$ and let $Q \in   X_i$ be the terminal endpoint of $\wt R_{E_q}$.     If $P = Q$ then $\ti R = \wt{R_{E_q}}$ and $\wt{R_{E_q}} \setminus \ti E_q \subset G_{r(j(i)-1)} $ by Corollary~\ref{high enough edges}.  In this case, $E_q = E^t$  so  $\sigma_0^{-1} =  R_{E^t}$ and    \pref{item: first term}  is satisfied.    Suppose then that $P \ne Q$.  The line $\ti L$  connecting $Q$ to $P$ is disjoint from $\Fix(\ti f_i)$.  Its highest edge splitting must be bi-infinite for otherwise the endpoints of its first or last term would be fixed.  It follows that  the height of  $\ti \rho$ equals the height of  $\wt{R_{E_q}}\setminus \ti E_q$ and so is at most $r(j(i)-1)$ by Corollary~\ref{high enough edges}.  It then follows that $E_q  =  E^t$ and $R_{E_t} = E^t w^{\pm\infty}$.  There is a conjugate $a'$ of $a$ such that $A_{a'}$ shares an endpoint with $\wt{R_{E_q}}$.    Lemma~\ref{line without fixed points} implies that $\ti L = A_{a'}$.   We conclude that $\ti \sigma_0^{-1}$ begins with $\ti E_q = \ti{E^t}$ and is otherwise contained in $A_{a'}$.  This completes the proof of   \pref{item: first term}.
   \endproof
   
   \begin{definition} \label{crossing number} If a path $\tau$ crosses $E^i$ and $\bar E^i$ a finite number of times then we define the {\em algebraic crossing number $c_i(\tau)$ of $\tau$ with $E^i$} to be the number of times that $\tau$ crosses $E^i$ minus the number of times that $\tau$ crosses $\bar E^i$.  For closed paths $\delta \subset G$, the formula $$\delta \mapsto c_i(\delta)\mod 3$$   defines a homomorphism       $$H_1(G;\Z/3) \mapsto \Z/3$$  and this homomorphism is nontrivial if and only if $E^i $ is nonseparating. For infinite paths of the form $\tau = w^\infty \delta w^\infty$, since $w$ does not cross $E^i$ we have $c_i(\tau) = c_i(\delta)$.
\end{definition}

   We sometimes use edge path notation in describing lines and rays.  In particular, a line or ray $\sigma$ {\em ends with $w^\infty$}  if the ray $R(w):= w^\infty$  is a 
terminal subray of $\sigma$ and  $\sigma$ {\em ends with $w^{-\infty}$}  if the ray $R(\bar w):= \bar w^\infty$  is a terminal subray of $\sigma$.  Analogously $\sigma$ {\em begins with $w^\infty  [resp. \ w^{- \infty}] $} if $\sigma^{-1}$ ends with $R(\bar w)$ [resp.$R(w)$]. 
   
   \begin{corollary}\label{index summary}   For all $1 \le i \le m-1$ and all $\sigma \in \Sigma_i$, \  $c_i(\sigma) = 0$ or $1$. 
\end{corollary}

\proof Let  $ \sigma = \cdot \ldots   \cdot \sigma_{-1}  \cdot   \sigma_0  \cdot   \sigma_1  \cdot  \ldots$ be the splitting given by Lemma~\ref{splitting at fixed points}. The corollary follows from
\begin{itemize}
\item $c_i(\sigma_s) = 0$ for each finite $\sigma_s$.
\item There is a last terminal term and its contribution to $c_i(\sigma)$   is $1$.  
\item If there is a first term then its  contribution to $c_i(\sigma)$   is $0$ or $-1$.  
\end{itemize}
each of which is an immediate consequence of Lemma~\ref{splitting at fixed points}.
\endproof

\noindent{\bf Verification that $\pi_a$ (as defined in $(*)$) is the identity:} \ \  To prove that $\pi_a$ is the identity, it suffices to show that $\pi_a(l) = l$ for all $l > 0$; the $l=0$ case then follows from the fact that $\pi_a$ is a permutation.     
  Suppose to the contrary that $l > 0$ and that  $\pi_a(k) =l$  for some $k\ne l$.

We claim that $\Sigma_k$ contains a line $\tau$ that begins with $ w^{\infty} $,   ends with $w^\infty$ and satisfies  $c_l(\tau) =  -1$. 
 If each line in $\Sigma_k$ is contained in $G \setminus E^l$  then the same is true for the realization in $G$ of each line in $\F(\Sigma_k)$    in contradiction to the fact (Lemma~\ref{splitting at fixed points}-\pref{item:last term infinite}) that every line in $\Sigma_l$ crosses $E^l$ and the fact (Lemma~\ref{Fsupp}-\pref{item:supp Sigma}) that $\F(\Sigma_l) = \F(\Sigma_k)$.  We may therefore choose a line $\mu \in \Sigma_k$ that crosses $E^l$ or $\bar E^l$.  By Remark~\ref{plus minus}, we may assume that $\mu$ ends with $w^\infty$.  Let $ \mu  = \cdot \ldots   \cdot \mu_{-1}  \cdot   \mu_0  \cdot   \mu_1  \cdot  \ldots$ be the splitting given by Lemma~\ref{splitting at fixed points}.   
The  last term in the splitting, if it exists,   does not cross $E_l$ or $\bar E_l$ by items \pref{item:last term infinite} and \pref{item:no last term} of Lemma~\ref{splitting at fixed points}. If some finite term $\mu_s$ crosses   $E^l$ or $\bar E^l$ then $\mu_s = E^lw^* \bar E^l$ by Lemma~\ref{splitting at fixed points}-\pref{item:finite term}.  In this case, $ \tau= w^\infty \bar E^l\mu_{s+1} \mu_{s+2}\ldots$ is contained in $\Sigma_k$ by Lemma~\ref{splitting at fixed points}-\pref{item:finite term}  and $c_l(\tau) = -1$.  If no finite term of the splitting crosses   $E^l$ or $\bar E^l$ then there must be a first term $\mu_b$ and  it must cross $E^l$ or $\bar E^l$   so   $\mu_b = w^{\pm \infty}\bar E^l$ by Lemma~\ref{splitting at fixed points}\pref{item: first term}.  In this case,  $ \tau= w^\infty \bar E^l\mu_{b+1} \mu_{b+2}\ldots$   satisfies $c_l(\tau) = -1$.   If $\tau \ne \mu$ then $\tau$ and $\mu$  have lifts  $\ti \tau$ and $\ti \mu$ with  terminal endpoint $a^+$     and with   initial endpoints   bounding $A_b$ some $b \in F_n$ satisfying $[b] =[w]$.  Since one of the endpoints of $A_b$ is contained in $X_k$ the other is also so $\tau \in \Sigma_k$ as desired.  This completes the proof of the claim.  

There is a closed path $\delta$ such that $\tau = w^\infty \delta w^\infty$.  The line $\tau' = h_\#(\tau)    \in \Sigma_l$ is obtained from $w^\infty h_\#(\delta) w^\infty$ by tightening.  No copies of $E_l$ or $\bar E_l$ are cancelled during the tightening process, so $\tau' = w^\infty \delta' w^\infty$ where $\delta'$ is a closed path satisfying $c_l(\delta') = c_l(h_\#(\delta))$.     Applying Definition~\ref{crossing number}, we have
$$c_l(\delta) =c_l(\tau) = -1 \qquad  \text{ and } \qquad c_l(h_\#(\delta)) = c_l(\delta') = c_l(\tau')$$  Since $\psi \in \IAThree$ and $h$ is a topological representative of $\psi$, 
 $$c_l(\tau') \text{ mod } 3 = c_l(h_\#(\delta))\text{ mod } 3 = c_l(\delta)\text{ mod } 3 \  = \ -1$$
         This contradiction to Corollary~\ref{index summary} completes the verification of ($*$) and hence the proof of Proposition~\ref{upg prop}.  \qed

\bibliographystyle{amsalpha} 
\bibliography{mosher} 

\end{document}